\documentclass[aps,preprint]{revtex4}

\usepackage[dvips]{graphicx}

\newcommand{\ch}{{\cal H}}
\newcommand{\ce}{{\cal E}}
\newcommand{\cpt}{PT}

\begin{document}

\date{\today}

\title{Exceptional points in open and PT symmetric systems}
\author{Hichem Eleuch$^{1}$\footnote{email: hichemeleuch@yahoo.fr} and 
Ingrid Rotter$^{2}$\footnote{email: rotter@pks.mpg.de}}

\address{
$^1$Department of Physics, McGill University, Montreal, Canada H3A 2T8}
\address{
$^2$Max Planck Institute for the Physics of Complex Systems,
D-01187 Dresden, Germany }

\begin{abstract}

Exceptional points (EPs) determine  the dynamics of open quantum
systems and cause also PT symmetry breaking in PT symmetric systems. From
a mathematical point of view, this is caused by the fact that the
phases of the wavefunctions (eigenfunctions of a non-Hermitian
Hamiltonian) relative to one another are not rigid when
an EP is approached. The  system is therefore able to align with the
environment to which it is coupled and, consequently, rigorous 
changes of the system properties may occur. 
We compare analytically as well as numerically the eigenvalues and
eigenfunctions of a $2\times 2$  matrix that is characteristic of either 
open quantum systems at high level density or of PT symmetric optical 
lattices. In both cases, the results show clearly the influence of the
environment onto the system in the neighborhood of EPs. Although the
systems are very different from one another, the eigenvalues and
eigenfunctions indicate the same characteristic features.

\end{abstract}

\pacs{\bf }
\maketitle

\section{Introduction}
\label{intr}

Starting with the paper \cite{bender1}, it
has been shown  that a wide class of PT symmetric non-Hermitian 
Hamilton operators provides entirely  real spectra. In the following
years this phenomenon is studied in 
many theoretical papers, see the review \cite{bender2} and the Special
Issue \cite{specjpa}. 

In order to realize complex PT symmetric structures, the formal
equivalence of the quantum mechanical Schr\"odinger equation to 
the optical wave equation in PT symmetric optical lattices
\cite{equivalence} can be exploited by involving symmetric index
guiding and an antisymmetric gain/loss profile. The  experimental
results \cite{mumbai1} have confirmed the expectations and have,
furthermore, demonstrated  the onset of passive PT symmetry 
breaking within the
context of optics. This phase transition was found to lead to a
loss-induced optical transparency in specially designed
pseudo-Hermitian potentials. In another experiment \cite{mumbai2}, 
the wave propagation in an active PT symmetric coupled waveguide system   
is studied. Both spontaneous PT symmetry breaking and power
oscillations violating left-right symmetry are observed. Moreover, the
relation of the relative phases of the eigenstates of the system to
their distance from the level crossing point is obtained. The phase
transition occurs when this point is approached. The  meaning of
these results for a new generation of integrated photonic devices is
discussed in \cite{kottos1}. Today we have 
many experimental and theoretical studies related to this topic.

On the other hand, non-Hermitian operators are known to describe open
quantum systems in a natural manner, see e.g. \cite{feshbach}. In
contrast to the original papers more than 50 years ago, 
statistical assumptions on the
system's states are not at all necessary today \cite{ro91} due to the
improved accuracy of the experimental as well as theoretical studies. 
In the present-day papers, the system is
assumed to be open due to the fact that it is embedded into the
continuum of scattering wavefunctions into which the states 
of the system can decay. This environment
exists always. It can be changed by means of external
forces, but cannot be deleted \cite{top}. 
The states of the system can decay due to their coupling to
the environment of scattering wavefunctions but cannot be formed out
of the continuum. Hence, loss  is nonvanishing usually,
while  gain is  zero. The complex eigenvalues of the non-Hermitian
Hamiltonian provide both the energy $E_i$ as well as the lifetime 
$\tau_i$ (inverse proportional to the decay width $\Gamma_i$) of the
eigenstate $i$.

Recent studies have shown the important role singular points in the
continuum  play for the dynamics of open quantum systems, see e.g. the
review \cite{top}. These
singular points are called usually exceptional points (EPs) after Kato
who studied their mathematical properties \cite{kato} many years ago. 
The relation of EPs  to PT symmetry breaking in optical systems
is considered already in the first papers  \cite{mumbai2,kottos1}.
Nevertheless, the relation between the 
dynamical properties of open quantum systems and
those of PT symmetric systems is not considered thoroughly up to now.

It is the aim of the present paper to compare directly the influence
of EPs onto the dynamics of open quantum systems with that onto 
PT symmetry breaking in PT symmetric systems. The comparison is
performed  on the basis of simple models with only two levels
coupled to one common channel. In both cases, 
the Hamiltonian is given by a  $2\times 2$ matrix in the form it is
used usually in the literature. We will follow here the representation
given for open quantum systems in \cite{top} and for PT symmetric
systems used in  \cite{jopt}.

In Sect. \ref{two1}, the non-Hermitian Hamiltonian of an open quantum
system is considered. The properties of its eigenvalues and eigenfunctions 
are sketched, above all in the neighborhood of one or more EPs.
In the following section \ref{two2}, two different non-Hermitian
operators that are used in the description of PT symmetric systems,  
are considered. The similarities and differences to the Hamiltonian
of an open quantum system are discussed on the basis of 
analytical studies (when possible) as well as by means of  numerical
results. The results  are summarized in the last section.

\section{Exceptional points in an open quantum system}
\label{two1}

In an open quantum system, the discrete states  described by a
Hermitian Hamiltonian $H^B$, are embedded into the continuum of scattering
wavefunctions, which exists always and can not be deleted. 
Due to this fact the discrete states turn into
resonance states the lifetime of which is usually finite. 
The Hamiltonian $\ch$ of the whole system consisting of the two subsystems, 
is non-Hermitian. Its eigenvalues are complex and
provide not only the energies  of the states but also their lifetimes
(being inverse proportional to the widths).
  
The Hamiltonian of an open quantum system reads \cite{top}
\begin{eqnarray}
\label{ham1}
\ch & = & H^B + V_{BC} G_C^{(+)} V_{CB} 
\end{eqnarray}
where $V_{BC}$ and $V_{CB}$ stand for the interaction between system
and environment and $  G_C^{(+)} $ is the Green function in the
environment. The so-called internal (first-order) interaction 
between two states $i$ and $j$ is
involved in $H^B$ while their external (second-order) interaction via the
common environment is described by the last term of (\ref{ham1}).

Generally, the coupling matrix elements of the external interaction consist of 
the principal value integral    
\begin{eqnarray}
{\rm Re}\; 
\langle \Phi_i^{B} | \ch |  \Phi_j^{B} \rangle 
 -  E_i^B \delta_{ij} =\frac{1}{2\pi} 
 {\cal P} \int_{\epsilon_c}^{\epsilon_{c}'} 
 {\rm d} E' \;  
\frac{\gamma_{ic}^0 \gamma_{jc}^0}{E-E'} 
\label{form11}
\end{eqnarray}
which is real, and the residuum
\begin{eqnarray}
{\rm Im}\; \langle \Phi_i^{B} | \ch |
  \Phi_j^{B} \rangle =
- \frac{1}{2}\; 
 \gamma_{ic}^0 \gamma_{jc}^0 
\label{form12}
\end{eqnarray}
which is imaginary \cite{top}. Here, the $\Phi_i^{B}$ and  $E_i^B$ are the
eigenfunctions and (discrete) eigenvalues, respectively, of
the Hermitian Hamiltonian $H^B$ which describes the states in the subspace of
discrete states without any interaction of the states via the
environment. The
$\gamma_{i c}^0  \equiv
\sqrt{2\pi}\, \langle \Phi_i^B| V | \xi^{E}_{c}
\rangle $
 are the (energy-dependent) coupling matrix elements  
between the discrete states $i$ of the system and the environment of
scattering wavefunctions $\xi_c^E$. The $\gamma_{k c}^0$ have to be
calculated for every state $i$ and for each channel $c$ 
(for details see \cite{top}). 
When $i=j$, (\ref{form11}) and (\ref{form12}) give the selfenergy of
the state $i$. 
The coupling matrix elements (\ref{form11}) and (\ref{form12})
(by adding $E_i^B \delta_{ij}$ in the first case) 
are often simulated by complex values  $\omega_{ij}$.

In order to study the interaction of two states via one common environment it
is convenient to start from two resonance states (instead of two
discrete states). 
Let us consider, as an example, the symmetric $2\times 2$ matrix 
\begin{eqnarray}
{\cal H}^{(2)} = 
\left( \begin{array}{cc}
\varepsilon_{1} \equiv e_1 + \frac{i}{2} \gamma_1  & ~~~~\omega_{12}   \\
\omega_{21} & ~~~~\varepsilon_{2} \equiv e_2 + \frac{i}{2} \gamma_2   \\
\end{array} \right) 
\label{form1}
\end{eqnarray}
the diagonal elements of which are the two complex eigenvalues 
$ \varepsilon_{i}~(i=1,2)$ of a non-Hermitian operator ${\cal H}^0$.
That means, the $e_i$ and  $\gamma_i \le 0$ denote the 
energies and widths, respectively, of the two states when
$\omega_{ij} =0$ (the index $c$ is ignored here for simplicity, $c=1$). 
The $\omega_{12}=\omega_{21}\equiv \omega$ stand for
the coupling of the two states via the common environment. The
selfenergy of the states is assumed to be included into the $\varepsilon_i$.

The two eigenvalues of ${\cal H}^{(2)}$ are
\begin{eqnarray}
\ce_{i,j} \equiv E_{i,j} + \frac{i}{2} \Gamma_{i,j} = 
 \frac{\varepsilon_1 + \varepsilon_2}{2} \pm Z ~; \quad \quad
Z \equiv \frac{1}{2} \sqrt{(\varepsilon_1 - \varepsilon_2)^2 + 4 \omega^2}
\label{int6}
\end{eqnarray}
where   $E_i$ and $\Gamma_i \le 0$ stand for the
energy and width, respectively, of the eigenstate $i$. 
Resonance states with nonvanishing widths $\Gamma_i$ 
repel each other in energy  according to the value of Re$(Z) $
while the widths bifurcate according to the value of Im$(Z)  $.
The two states cross when $Z=0$. This crossing point is 
an EP according to the definition of Kato \cite{kato}.
Here, the two eigenvalues  coalesce, $\ce_{1}=\ce_{2}$.

According to (\ref{int6}), two interacting discrete states (with
$\gamma_1 = \gamma_2 =   0$) avoid always crossing since 
$\omega \equiv \omega_0$ and 
$\varepsilon_1 - \varepsilon_2$ are real in this case and the 
condition $Z=0$ can not be fulfilled,
\begin{eqnarray}
(e_1 - e_2)^2 +4\, \omega_0^2 &>& 0.  
\label{int6a}
\end{eqnarray}
In this case, the EP can be found only by
analytical continuation into the continuum. This situation is known as
avoided crossing of discrete states.
It holds also for narrow resonance states if $Z=0$ cannot be
fulfilled due to the small widths of the two states. 
The physical meaning of this result is very well known since many
years. The avoided crossing of two discrete states at a certain critical
parameter value \cite{landau} means that the
two states are exchanged at this  point, including their 
populations ({\it  population transfer}).  
 
When $\omega = i ~\omega_0 $ is imaginary, 
\begin{eqnarray}
Z = \frac{1}{2} \sqrt{(e_1-e_2)^2 - \frac{1}{4} (\gamma_1-\gamma_2)^2 
+i(e_1-e_2)(\gamma_1-\gamma_2) - 4\omega_0^2}
\label{int6i}
\end{eqnarray}
is complex. The condition  $Z= 0$ can be fulfilled only when $
(e_1-e_2)^2 - \frac{1}{4} (\gamma_1-\gamma_2)^2 = 4\omega_0^2$
and $(e_1-e_2)(\gamma_1-\gamma_2) =0$, i.e. 
when $\gamma_1 = \gamma_2$ (or when $e_1=e_2$). In this case, it follows 
\begin{eqnarray}
(e_1 - e_2)^2 -4\, \omega_0^2 &= &0 
~~\rightarrow ~~e_1 - e_2 =\pm \, 2\, \omega_0 
\label{int6b}
\end{eqnarray}
and two EPs appear.  It holds further
\begin{eqnarray}
\label{int6c}
(e_1 - e_2)^2 >4\, \omega_0^2 &\rightarrow& ~Z ~\in ~\Re \\
\label{int6d}
(e_1 - e_2)^2 <4\, \omega_0^2 &\rightarrow&  ~Z ~\in ~\Im 
\end{eqnarray}
independent of the parameter dependence of the $e_i$.
In the first case, the eigenvalues ${\cal E}_i = E_i+i/2\, \Gamma_i$ 
differ from the original values 
$\varepsilon_i = e_i + i/2~\gamma_i$ by a contribution to
the energies and in the second case by a contribution to the widths. 
The width bifurcation starts in the very neighborhood of one of the EPs and
becomes maximum in the middle between the two EPs.
This happens at the crossing point $e_1 = e_2$ where 
$\Delta \Gamma/2 \equiv |\Gamma_1/2 - \Gamma_2/2| = 4\, \omega_0$.
A similar situation appears when $\gamma_1 \approx \gamma_2$ as 
results of numerical calculations  show.
The physical meaning of this result is completely different from that
discussed above for discrete and narrow resonance states.
It means that {\it different time scales} appear in the system without any
enhancement of the coupling strength to the continuum (for details see
\cite{fdp1}).

The cross section can be calculated by means of the $S$ matrix 
$\sigma (E) \propto |1-S(E)|^2$.  A unitary representation of the 
$S$ matrix in the case of two nearby resonance states coupled to one 
common continuum of scattering wavefunctions reads \cite{top} 
\begin{eqnarray}
\label{sm1}
S = \frac{(E-E_1+\frac{i}{2}\Gamma_1)~(E-E_2+\frac{i}{2}\Gamma_2)}{(E-E_1-
\frac{i}{2}\Gamma_1)~(E-E_2-\frac{i}{2}\Gamma_2)} \; .
\end{eqnarray}
In this expression, the influence of an EP onto the cross section 
is contained in the eigenvalues 
${\cal{E}}_i = E_i + i/2~\Gamma_i$ of $\ch^{(2)}$.
Reliable results can be obtained therefore  also when 
an EP is approached and the $S$ matrix has a double pole. Here, the
line shape of the two overlapping resonances is described by
\begin{eqnarray}
\label{sm2}
S = 1+2i\frac{\Gamma_d}{E-E_d-\frac{i}{2}\Gamma_d}-
\frac{\Gamma_d^2}{(E-E_d-\frac{i}{2}\Gamma_d)^2}
\end{eqnarray}
where  $E_1=E_2\equiv E_d$  and $\Gamma_1=\Gamma_2\equiv \Gamma_d$.
It deviates  from the Breit-Wigner line shape of an isolated resonance 
due to  interferences between the two resonances. The first term of 
(\ref{sm2}) is linear (with the factor $2$ in front) while the second one is  
quadratic. As a result, two peaks with asymmetric line shape
appear in the cross section (for a numerical example see Fig. 9 in 
\cite{mudiisro}).

The eigenfunctions of the non-Hermitian $\ch^{(2)}$ are biorthogonal 
and can be normalized according to
\begin{eqnarray}
\langle \Phi_i^*|\Phi_j\rangle = \delta_{ij} 
\label{int3}
\end{eqnarray}
although $\langle \Phi_i^*|\Phi_j\rangle  $ is a complex number
(for details see sections 2.2 and 2.3 of \cite{top}). The
normalization (\ref{int3}) allows to describe the smooth transition from the
regime with orthogonal eigenfunctions to that with biorthogonal  
eigenfunctions (see below). It follows 
\begin{eqnarray}
 \langle\Phi_i|\Phi_i\rangle & = & 
{\rm Re}~(\langle\Phi_i|\Phi_i\rangle) ~; \quad
A_i \equiv \langle\Phi_i|\Phi_i\rangle \ge 1
\label{int4} 
\end{eqnarray}
and 
\begin{eqnarray}
\langle\Phi_i|\Phi_{j\ne i}\rangle & = &
i ~{\rm Im}~(\langle\Phi_i|\Phi_{j \ne i}\rangle) =
-\langle\Phi_{j \ne i}|\Phi_i\rangle 
\nonumber  \\
&& |B_i^j|  \equiv 
|\langle \Phi_i | \Phi_{j \ne i}| ~\ge ~0  \; .
\label{int5}
\end{eqnarray}
At an EP $A_i \to \infty$ and $|B_i^j| \to \infty$.
The $\ce_i$ and $\Phi_i$ contain global features that are 
caused by many-body forces  induced by the coupling
$\omega_{ik}$ of the states $i$ and $k\ne i$ via the environment.
They contain moreover the self-energy  of the states $i$
due to their coupling to the environment. 

At the EP,
the eigenfunctions $\Phi_i^{\rm cr}$ of ${\cal H}^{(2)}$
of the two crossing states are linearly dependent from one another, 
\begin{eqnarray}
\Phi_1^{\rm cr} \to ~\pm ~i~\Phi_2^{\rm cr} \; ;
\quad \qquad \Phi_2^{\rm cr} \to
~\mp ~i~\Phi_1^{\rm cr}   
\label{eif5}
\end{eqnarray}  
according to analytical  as well as numerical and experimental
studies, see  Appendix of \cite{fdp1} and 
section 2.5 of \cite{top}.
This means, that the wavefunction $\Phi_1$ of the state $1$ jumps, 
at the EP, via the wavefunction ~$\Phi_1\pm \, i\, \Phi_2$  
of a chiral state to   ~$\pm\, i\, \Phi_2$
\cite{comment}. 

The Schr\"odinger equation with the non-Hermitian operator 
${\cal H}^{(2)}$ is equivalent to a Schr\"odinger equation with 
${\cal H}^0$ and source term \cite{ro01}
\begin{eqnarray}
\label{form1a}
({\cal H}^0 - \varepsilon_i) ~| \Phi_i \rangle  = -
\left(
\begin{array}{cc}
0 & \omega_{ij} \\
\omega_{ji} & 0
\end{array} \right) |\Phi_j \rangle \equiv W  |\Phi_j \rangle\; . 
\end{eqnarray}
Due to the source term, two states are coupled via the 
common environment of scattering wavefunctions into which the system 
is embedded,  $\omega_{ij}=\omega_{ji}\equiv\omega$.

The Schr\"odinger equation (\ref{form1a}) with source term can be
rewritten in the following manner \cite{ro01},
\begin{eqnarray}
\label{form2a}
({\cal H}^0  - \varepsilon_i) ~| \Phi_i \rangle  = 
\sum_{k=1,2} \langle
\Phi_k|W|\Phi_i\rangle \sum_{m=1,2} \langle \Phi_k |\Phi_m\rangle 
|\Phi_m\rangle \; . 
\end{eqnarray}
According to the biorthogonality  relations
(\ref{int4}) and (\ref{int5}) of the eigenfunctions of ${\cal H}^{(2)}$,  
(\ref{form2a}) is a nonlinear equation.  
Most important part of the nonlinear contributions is contained in 
\begin{eqnarray}
\label{form3a}
({\cal H}^0  - \varepsilon_n) ~| \Phi_n \rangle =
\langle \Phi_n|W|\Phi_n\rangle ~|\Phi_n|^2 ~|\Phi_n\rangle \; .  
\end{eqnarray}
The nonlinear source term vanishes far from an EP due to
$\langle \Phi_k|\Phi_{k }\rangle    \to  1 $ and
$\langle \Phi_k|\Phi_{l\ne k }\rangle = - 
\langle \Phi_{l \ne k  }|\Phi_{k}\rangle \to  0 $ according to 
(\ref{int3}) to  (\ref{int5}).
Thus, the Schr\"odinger equation with source term is 
linear far from an EP, as usually assumed. It is however nonlinear
in the neighborhood of an EP.

It is meaningful to represent
the  eigenfunctions $\Phi_i$ of ${\cal H}^{(2)}$  in the
set of basic wavefunctions $\Phi_i^0$ of ${\cal H}^0$
\begin{eqnarray}
\Phi_i=\sum_{j=1}^N b_{ij} \Phi_j^0 ~~ ;
\quad \quad b_{ij} = |b_{ij}| e^{i\theta_{ij}}
\; .
\label{int20}
\end{eqnarray}
Also the $b_{ij}$ are normalized  according to the biorthogonality
relations  of the wavefunctions $\{\Phi_i\}$. The angle $\theta_{ij}$
can be determined from
${\rm tg}(\theta_{ij}) = {\rm Im}(b_{ij}) / {\rm Re}(b_{ij})$ .

From (\ref{int3}) and (\ref{eif5}) follows\,:
\begin{verse}
(i) When  two levels are distant from one another,  their eigenfunctions
 are (almost) orthogonal,  
$\langle \Phi_k^* | \Phi_k \rangle   \approx
\langle \Phi_k | \Phi_k \rangle  = A_k \approx 1 $.\\
(ii) When  two levels cross at the EP, their eigenfunctions are linearly
dependent according to (\ref{eif5}) and 
$\langle \Phi_k | \Phi_k \rangle \equiv A_k \to \infty $.\\
\end{verse}
These two relations show that the phases of the two eigenfunctions
relative to one another change when the crossing point is approached. 
This can be expressed quantitatively by defining the {\it phase
  rigidity} $r_k$ of the eigenfunction $\Phi_k$,
\begin{eqnarray}
r_k ~\equiv ~\frac{\langle \Phi_k^* | \Phi_k \rangle}{\langle \Phi_k 
| \Phi_k \rangle} ~= ~A_k^{-1} \; . 
\label{eif11}
\end{eqnarray}
It holds $1 ~\ge ~r_k ~\ge ~0 $.  
The  non-rigidity $r_k$ of the phases of the eigenfunctions of $\ch^{(2)}$ 
follows also from the fact that $\langle\Phi_k^*|\Phi_k\rangle$
is a complex number (in difference to the norm
$\langle\Phi_k|\Phi_k\rangle$ which is a real number) 
such that the normalization condition
(\ref{int3}) can be fulfilled only by the additional postulation 
Im$\langle\Phi_k^*|\Phi_k\rangle =0$ (what corresponds to a rotation, 
generally). 

When $r_k<1$, an analytical expression for the eigenfunctions as a
function of a certain control parameter  can, generally, not be
obtained. The  non-rigidity $r_k<1$ of the phases of the eigenfunctions
of $\ch^{(2)}$ in the neighborhood of EPs is the most important
difference between the  non-Hermitian quantum physics and the Hermitian
one. Mathematically, it causes nonlinear effects in quantum systems 
in a natural manner, as shown above.
Physically, it allows the alignment of one of the states of the
system to the common environment \cite{top}.

Results of numerical calculations are given in, e.g.,  \cite{elro2}.
The mixing coefficients $b_{ij}$  (defined in (\ref{int20}))
of the wavefunctions of the two  states due to their avoided crossing 
are simulated by assuming a  Gaussian distribution 
for the coupling coefficients $\omega_{i\ne j} = \omega ~e^{-(e_i  -e_j)^2} $ 
(for real $\omega$, the results of the simulation agree with the
results \cite{ro01} of exact calculations). In \cite{elro2}, 
results of different calculations are shown for illustration. Here, 
the coupling coefficients $\omega$ are assumed to be either real or 
complex or imaginary according to the different possibilities provided by
(\ref{form11}) and (\ref{form12}).

The main difference of the eigenvalue trajectories with
real to those with imaginary coupling coefficients $\omega$ is
related to the relations (\ref{int6a}) to (\ref{int6d}) obtained
analytically.
For $\gamma_1 \ne \gamma_2$ and real, complex or even imaginary 
$\omega$, the results show one EP when the condition $Z=0$
is fulfilled. This EP is isolated from
other EPs, generally, when the level density is low. In the case of
$\gamma_1 \approx \gamma_2$ and imaginary $\omega$  
however, two related EPs appear, see Fig. \ref{fig1} right panel.
Between these two EPs, the widths $\Gamma_i$ bifurcate
(Fig. \ref{fig1}.d) while the energies $E_i$ do not change
(Fig. \ref{fig1}.b).
It is interesting to see that  width bifurcation occurs {\it between}
the two EPs, according to (\ref{int6b}) and (\ref{int6d}),
{\it without} any enhancement of the coupling strength to the environment.
Beyond the two EPs, the eigenvalues approach the original values.

In a finite neighborhood of the point at which the two eigenvalue 
trajectories cross, the eigenfunctions 
are mixed and  $|b_{ij}|\to \infty$ when approaching the EP 
(Fig. \ref{fig1}.f).  The phases of {\it all} components of the
eigenfunctions jump at the EP either by $-\pi /4$ or by
$+\pi /4$ \cite{elro4}. That means the phases of {\it both} eigenfunctions
jump in the same direction by the same amount. Thus, there is a phase
jump of $-\pi /2$ (or $+\pi /2$) when one  of the eigenfunctions passes into   
the other one at the EP. This result is in agreement with (\ref{eif5}).
It holds true for real as well as for imaginary  $\omega$.

\section{Exceptional points in PT symmetric systems}
\label{two2}

As has been shown in \cite{equivalence}, 
the optical wave equation in $\cpt$ symmetric optical lattices 
is formally equivalent to a quantum mechanical Schr\"odinger equation. 
Complex $\cpt$ symmetric structures can be realized by involving 
symmetric index guiding and an antisymmetric gain/loss profile.

\begin{figure}[ht]
\begin{center}
\includegraphics[width=13cm,height=13cm]{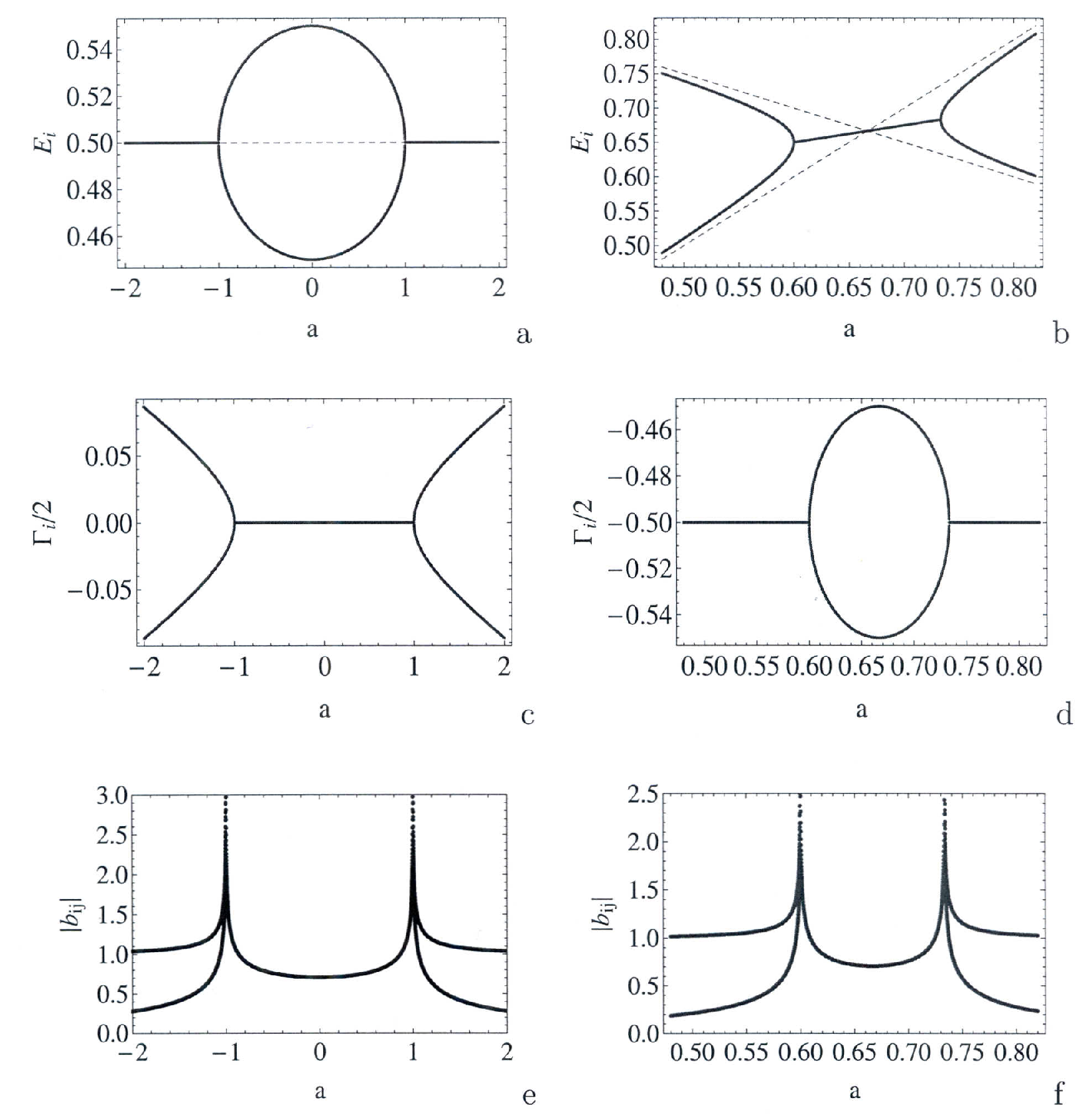}
\end{center}
\caption{\footnotesize
Energies $E_i$,  widths $\Gamma_i/2$ and wavefunctions $|b_{ij}|$
of $N=2$ states coupled to
$K=1$ channel as  function of $a$ of a PT symmetric system with 
Hamiltonian (\ref{ham2o1}) (left panel) 
and of an open quantum system with Hamiltonian (\ref{form1}) (right panel).
Parameters left panel: $e=0.5; ~\gamma_1=-\gamma_2=0.05~a; ~w =0.05$; 
~right panel: $e_1=1-0.5~a; ~e_2=a$; ~$\gamma_1/2= \gamma_2/2=0.5;
~\omega = 0.05\; i$. The dashed lines in (a,b) show $e_i(a)$.
}
\label{fig1}
\end{figure}

The main difference of
these optical systems to open quantum systems consists in the asymmetry of
gain and loss in the first case while the states 
of an open quantum system can only decay (Im$(\varepsilon_{1,2}) < 0$
and Im$({\cal E}_{1,2}) <0$
for all states).  Thus, the  modes involved in the non-Hermitian 
Hamiltonian in optics appear in complex conjugate pairs 
while this is not the case in an open quantum system.
As a consequence, the Hamiltonian for $\cpt$ symmetric structures in 
optical lattices may have real eigenvalues in a large parameter range.
The $2\times 2$ non-Hermitian Hamiltonian may be written as 
\cite{equivalence,jopt}
\begin{eqnarray}
\label{ham2o1}
\ch_{PT} =\left(
\begin{array}{cc}
~e-i\frac{\gamma}{2}~ &  w \\ w^* & ~e+i\frac{\gamma}{2}~
\end{array}
\right), 
\end{eqnarray}
where $ e $ stands for the energy  of the two modes, 
$\pm \gamma$ describes  gain and loss, respectively,
and the coupling coefficients $w $ stand for the
coupling of the two modes via the lattice. 
When the $\cpt$ symmetric optical lattices are studied with vanishing 
gain, the Hamiltonian reads 
\begin{eqnarray}
\label{ham2o2}
\ch_{PT}' =\left(
\begin{array}{cc}
~~e-i\frac{\gamma}{2}~~ &  w \\ w^* & e
\end{array}
\right) \; .
\end{eqnarray}
In realistic systems, $ w$ in (\ref{ham2o1}) and (\ref{ham2o2}) is
mostly  real (or almost real)  \cite{kottos2}.

The eigenvalues of the Hamiltonian (\ref{ham2o1}) differ from 
(\ref{int6}),
\begin{eqnarray}
\label{eig2o1}
\ce^{PT}_\pm &=& e \pm \frac{1}{2}\sqrt{4|w|^2 - \gamma^2} \equiv e
\pm  Z_{PT} \, .
\end{eqnarray}
A similar expression is derived in \cite{mumbai1}.
Since $e$ and $\gamma$ are real, the $\ce^{PT}_\pm$ are real when
$4|w|^2 > \gamma^2$. Under this condition,
the two levels repel each other in energy what is
characteristic of discrete interacting states. The level repulsion
decreases with increasing $\gamma$ (when the interaction $w$ is
fixed). When $4|w|^2 = \gamma^2 $  the two states cross.
Here,  $\ce^{PT}_\pm = e$ and $\gamma = \pm \sqrt{4|w|^2}$.
With further increasing $\gamma$ and $4|w|^2 < \gamma^2$
($w$ fixed for illustration), width bifurcation (PT symmetry breaking)
occurs and $\ce^{PT}_\pm = e \pm \frac{i}{2}\sqrt{\gamma^2 - 4|w|^2}$.

These relations are in
accordance with (\ref{int6b})  to (\ref{int6d}) for open quantum systems.
Two EPs exist according to  
\begin{eqnarray} 
\label{pt2}
4 |w|^2=(\pm \gamma)^2 \; . 
\end{eqnarray}
Further
\begin{eqnarray}
\label{pt3}
\gamma^2 <4\, |w|^2 &\rightarrow& ~Z_{PT} ~\in ~\Re \\
\label{pt4}
\gamma^2 >4\, |w|^2 &\rightarrow&  ~Z_{PT} ~\in ~\Im 
\end{eqnarray}
independent of the parameter dependence $\gamma(a)$ and of the ratio
Re$(w)$/Im$(w)$.  

In the case of the Hamiltonian (\ref{ham2o2}), the eigenvalues read
\begin{eqnarray}
\label{eig2o2}
\ce^{'PT}_\pm = e - i~\frac{\gamma}{4}
\pm \frac{1}{2}\sqrt{4|w|^2 - \frac{\gamma^2}{4}} \equiv e -
i~\frac{\gamma}{4} \pm Z_{PT}'.
\end{eqnarray}
We have level repulsion as long as 
$4|w|^2 > \frac{\gamma^2}{4}$. While  level repulsion 
decreases with increasing $\gamma$,
loss increases with increasing $\gamma$.
At the crossing point, $\ce^{'PT}_\pm = e - i~\frac{\gamma}{4}$. 
With further increasing $\gamma$ and $ 4|w|^2 \ll \frac{\gamma^2}{4}$ 
\begin{eqnarray}
\label{eig2o3}
\ce^{'PT}_\pm \to  e - i~\frac{\gamma}{4} \pm i~\frac{\gamma}{4}
= \Bigg\{    \begin{array}{c}
e \qquad \\
e - i~\frac{\gamma}{2}.
\end{array} 
\end{eqnarray}
The two modes (\ref{eig2o3}) behave differently. While the loss in one
of them is large, it is almost zero in the other one. Thus, only one
of the modes effectively survives. Equation (\ref{eig2o3})
corresponds to high transparency at large $\gamma$.

Further, two EPs exist according to 
\begin{eqnarray} 
\label{pt6}
4 |w|^2=(\pm \gamma /2)^2 
\end{eqnarray}
and
\begin{eqnarray}
\label{pt7}
\gamma^2/4 <4\, |w|^2 &\rightarrow& ~Z_{PT}^{'} ~\in ~\Re \\
\label{pt8}
\gamma^2/4 >4\, |w|^2 &\rightarrow&  ~Z_{PT}^{'} ~\in ~\Im \; . 
\end{eqnarray}
In analogy  to (\ref{pt2}) up to (\ref{pt4})), these relations 
are independent of the parameter dependence of $\gamma$ and of the ratio
Re$(w)$/Im$(w)$. 

Thus, the difference between the eigenvalues $\ce_i$ of $\ch^{(2)}$ of an  
open quantum system and the eigenvalues 
of the Hamiltonian of a 
PT symmetric system consists, above all, in the fact that  the $\ce_i$
depend on the ratio Re$(\omega)$/Im$(\omega)$ while the $\ce^{PT}_\pm$
and $\ce^{'PT}_\pm$ are independent of Re$(w)$/Im$(w)$. There exist
however similarities between the two cases.

Interesting is the comparison of the eigenvalues $\ce_i$ of $\ch^{(2)}$
obtained for imaginary non-diagonal matrix elements $\omega$, with the 
eigenvalues of (\ref{ham2o1}) (or (\ref{ham2o2})) obtained for real $w$.
In both cases, there are two EPs, see Fig. \ref{fig1}. 
In the first case (right panel), 
the energies $E_i$ of both states are equal and the widths $\Gamma_i$
bifurcate between the two EPs. This situation is
characteristic of an open quantum system at high level density
with complex (almost imaginary) $\omega$, see Eqs.
(\ref{int6b}) to (\ref{int6d}). In the second case (left panel) however 
the difference $|E_1 - E_2|$ of the energies increases (level
repulsion) first and decreases then again
while the  widths $\Gamma_i$ of both states vanish  in
the parameter range between the two EPs
in accordance with the analytical results (\ref{pt2}) to (\ref{pt4}). 
Between the two EPs, level repulsion causes
the two levels to be distant from one another and $w$
is expected to be (almost) real.  This result agrees qualitatively with
(\ref{form11}) 
and (\ref{form12}). Similar results are obtained for the eigenvalues of
(\ref{ham2o2}). The only difference to those of (\ref{ham2o1})
is that the $\Gamma_i$ do not vanish but
decrease between the two EPs with increasing $a$ in this case.

According to Figs. \ref{fig1}.a-d, the role of energy and width is formally
exchanged when the eigenvalues of the Hamiltonian (\ref{form1}) are
compared with those of (\ref{ham2o1}) (or (\ref{ham2o2})).
In any case, the eigenvalues are influenced strongly by the EPs.

Also the eigenfunctions of the Hamiltonian (\ref{form1}) of an open quantum
system (with imaginary $\omega$)
and those of the Hamiltonians (\ref{ham2o1}) and (\ref{ham2o2})
of a PT symmetric system (with real $w$) show similar features.
The eigenfunctions  $\Phi_i^{PT}$ of $\ch_{PT}$
(and $\Phi_i^{'PT}$ of $\ch_{PT}'$) are biorthogonal with all the
consequences discussed in Sect. \ref{two1}. In contrast to the
eigenvalues, they are dependent on the ratio  Re$(\omega)$/Im$(\omega)$.

The eigenfunctions can be represented  in a
set of basic wavefunctions in full analogy to the representation of
the eigenfunctions $\Phi_i$ of $\ch^{(2)}$ in (\ref{int20}).
They contain valuable information on the mixing of the wavefunctions
under the influence of the non-diagonal coupling matrix elements $w$
and $w^*$ in (\ref{ham2o1}) and (\ref{ham2o2}), respectively, and its relation
to EPs.  Due to the level repulsion occurring between the two EPs, 
the coupling coefficients $w$ can be considered to be (almost) real
in realistic cases. The phases of the
eigenmodes of the non-Hermitian Hamiltonians (\ref{ham2o1}) and (\ref{ham2o2})
are not  rigid, generally, in approaching an EP, and  
spectroscopic redistribution processes occur in the system  
under the influence of the environment (lattice). 
As in the case of open quantum systems,
the phase rigidity $r_k$ can be defined according to
(\ref{eif11}). It varies between 1 and 0 and is a quantitative measure
for the skewness of the modes when the crossing point is approached.

In Figs. \ref{fig1}.e and f,  the
eigenfunctions of the Hamiltonian (\ref{ham2o1})
(calculated with real $w$)  
are compared to those of the Hamiltonian (\ref{form1}) (calculated with 
imaginary $\omega$). They  show the same characteristic features.  
As can be seen from Fig. \ref{fig1}.e, PT symmetry breaking 
is accompanied by a mixing of the eigenfunctions {\it in a finite
  neighborhood} of the EPs in PT symmetric
systems. This result is  in complete analogy to the results shown 
in Fig. \ref{fig1}.f
for open quantum systems where a hint to width bifurcation
can be seen in the mixing of the eigenfunctions 
{\it around} these points. Also the phases of the
eigenfunctions jump in both cases by $\pi /4$ at the EPs (not shown here).
In the parameter region between the two EPs, the eigenfunctions
are completely mixed (1:1) in both cases while they
are unmixed far beyond the EPs, see Figs. \ref{fig1}.e and f.

\section{Discussion of the results}
\label{concl}

On the basis of $2\times 2$ models, we have compared the influence of
an EP onto the dynamics of an open quantum system with its influence onto PT
symmetry breaking  in a PT symmetric system. In the first case 
the coupling of the two states via the environment is symmetric
($\omega_{12} = \omega_{21} \equiv \omega$). In the second case however,
the formal equivalence of the optical wave equation in PT symmetric 
optical lattices with a quantum mechanical Schr\"odinger equation 
causes  the two nondiagonal matrix elements to be complex conjugate 
($w_{21} = w_{12}^*$). The eigenvalues depend in the first case on the
ratio Re${(\omega)}$/Im${(\omega)}$ while they are independent of 
Re${(w)}$/Im${(w)}$ in the second case.  The eigenfunctions are
sensitive to Re${(\omega)}$/Im${(\omega)}$  and Re${(w)}$/Im${(w)}$,
respectively, in both cases. 

The EPs cause nonlinear effects in their neighborhood which determine
the evolution  of open as well as of PT symmetric
systems. Most important for the dynamics of an open quantum system is
the regime at high level density where the coupling coefficients are
(almost) imaginary. Here, two EPs appear when the decay widths
$\gamma_i$ of both states are (almost) the same. Approaching the EPs,  width
bifurcation starts and ends, respectively, while beyond the EPs the
widths of both states are equal (or similar) to one another. The
energies of the two states show an opposite behavior: it is 
$E_i = E_2$ (or  $E_i\approx E_2$) in the parameter range between the
two EPs while the states repel each other in energy beyond the EPs.   
The width bifurcation related to the two EPs becomes relevant for the
dynamics of an open quantum system at high level density. Here,
short-lived and long-lived states are formed which are related to different 
time scales of the system (for details see \cite{fdp1}). 

Two EPs appear also in a PT symmetric system, and PT symmetry breaking is
directly related to them. From a mathematical point of view however, 
energy and time are exchanged in comparison with the corresponding
values  in an open quantum system. That means the widths of both states
are equal and vanish in the case of the Hamiltonian (\ref{ham2o1})
with gain and loss in the whole parameter range between the two EPs.
In this parameter range, the eigenvalues are real and, furthermore,
level repulsion prohibits a small energy distance
between the two levels. Therefore the non-diagonal coupling
matrix elements $w$ are (almost) real, Re$(w) \gg {\rm Im}(w)$.  

The eigenfunctions of the different  $2\times 2$ models 
considered in the present paper, show very
clearly that the spectroscopic redistribution
inside the system is caused by the EPs, indeed.
However, it shows up  in all cases in a finite neighborhood around
them. Here the rigidity of the phases of the two eigenfunctions
relative to one another is reduced ($r_i<1$) and an alignment of one
of the states to the environment is possible. In the parameter range   
between the two EPs, the wavefunctions are completely mixed (1:1) as
can be seen from the numerical results shown in Fig. \ref{fig1}.

Summing up the discussion we state the following. The results obtained
by studying PT symmetric optical lattices as well as those received
from an investigation of open quantum systems 
show the characteristic features of non-Hermitian quantum physics. 
They prove environmentally induced effects that cannot be described 
convincingly in conventional Hermitian quantum physics. 
Due to the reduced phase
rigidity around an EP, the system is able to align (at least partly)
with the environment. This can be seen from PT symmetry breaking
occurring in one of the considered systems as well as
from  the dynamical phase transition taking place at high level 
density in the other system.


\begin{thebibliography}{99}
\bibitem{comment} In studies by other researchers, the factor $i$ in 
(\ref{eif5}) does not appear. This difference is
discussed in detail and compared with experimental data in the
Appendix of \cite{fdp1} and in section 2.5 of \cite{top}
\end{thebibliography}
\end{document}